\documentclass[11pt,a4paper]{article}
\usepackage[T1]{fontenc}
\usepackage[utf8]{inputenc}
\usepackage{lmodern}
\usepackage[margin=1in]{geometry}
\usepackage{amsmath,amssymb,amsthm,color}
\usepackage{microtype}
\usepackage[hidelinks]{hyperref}
\hypersetup{
	pdftitle={Every graph contains an odd induced subgraph with linearly many edges},
	pdfauthor={Jiangdong Ai, Qiwen Guo, Gregory Gutin and Yacong Zhou},
	pdfkeywords={odd induced subgraph, even induced subgraph, local complementation, Schur complement, entropy}
}

\theoremstyle{plain}
\newtheorem{theorem}{Theorem}[section]
\newtheorem{lemma}[theorem]{Lemma}

\newtheorem{corollary}[theorem]{Corollary}
\newtheorem{question}[theorem]{Question}

\newcommand{\YZ}[1]{{#1}}

\newcommand{\F}{\mathbb{F}_2}
\newcommand{\A}{\mathsf{A}}
\newcommand{\one}{\mathbf{1}}
\newcommand{\cM}{\mathcal{M}}
\newcommand{\EE}{\mathbb{E}}
\DeclareMathOperator{\diag}{diag}
\DeclareMathOperator{\Diag}{Diag}
\DeclareMathOperator{\supp}{supp}
\DeclareMathOperator{\rank}{rank}

\title{Every graph contains an odd induced subgraph with linearly many edges}

\author{Jiangdong Ai\thanks{School of Mathematical Sciences and LPMC, Nankai University, Tianjin 300071, China. {\tt jd@nankai.edu.cn}. Funded by the National Natural Science Foundation of China (No.12522117, No.12401456), the Natural Science Foundation of Tianjin (No.24JCQNJC01960) and Fundamental and Interdisciplinary Disciplines Breakthrough Plan of the Ministry of Education of China (JYB2025XDXM207).}
	\and Qiwen Guo\thanks{Department of Computing, Security and Mathematics, Royal Holloway, University of London, Egham, Surrey TW20 0EX, UK. {\tt gqwmath@163.com}.}
	\and Gregory Gutin\thanks{Corresponding author. Department of Computing, Security and Mathematics, Royal Holloway, University of London, Egham, Surrey TW20 0EX, UK. {\tt g.gutin@rhul.ac.uk}.}
	\and Yacong Zhou\thanks{Shenzhen Institutes of Advanced Technology, Chinese Academy of Sciences, Shenzhen, China. {\tt yacong.zhou96@gmail.com}.}
}

\date{}

\begin{document}
	
	\maketitle
	
	\begin{abstract}
		A graph is odd (even, respectively) if all its vertices have odd (even, respectively) degree. 
		For a graph $G$ without isolated vertices, let $f_o(G)$ be the maximum number of vertices in an odd 
		induced subgraph of $G$. Solving a long-standing conjecture, Ferber and Krivelevich proved in 2022 that there is a constant $c>0$ such that $f_o(G)\ge cn$ for every graph $G$ with $n$ non-isolated vertices.
		For a graph $G$, let $f_o^*(G)$ be the maximum number of edges in an odd induced subgraph of $G$. Scott asked in 2001 whether there is a constant $c>0$ such that $f_o^*(G)\ge cm$ for every graph $G$ with $m$ edges. We answer this question in the affirmative. We prove that if $G$ has $m$ edges and $n$ non-isolated vertices, then $m\le \YZ{1928}f_o^*(G)+\YZ{88}n$, and consequently $f_o^*(G)\ge m/\YZ{3000}$. We also prove that every graph with $n$ vertices and $m$ edges has an even induced subgraph with at least $\YZ{(m-11n)/120}$ edges; a term linear in $n$ cannot be avoided here, since an even induced subgraph of a forest has no edges. 
		
		
	\end{abstract}
	
	\section{Introduction}\label{sec:intro}
	
	All graphs in this paper are finite and simple. A graph is \emph{odd} if every vertex has odd degree, and \emph{even} if every vertex has even degree. A classical theorem of Gallai (see~\cite[Problem~5.17]{Lovasz}) states that the vertex set of every graph can be partitioned into two sets, each of which induces an even subgraph. In particular, every graph on $n$ vertices has an even induced subgraph on at least $n/2$ vertices. For odd induced subgraphs the situation is different, since an isolated vertex belongs to no odd subgraph. For a graph $G$, let $f_o(G)$ be the maximum order of an odd induced subgraph of $G$. It was a long-standing conjecture, stated by Caro~\cite{Caro}, that there is a constant $c>0$ such that $f_o(G)\ge cn$ for every graph $G$ on $n$ vertices without isolated vertices. Answering a question of Alon, Caro~\cite{Caro} proved that $f_o(G)=\Omega(\sqrt n)$, and Scott~\cite{Scott92} improved this to $f_o(G)=\Omega(n/\log n)$. The conjecture was proved by Ferber and Krivelevich~\cite{FK} with $c=10^{-4}$. The constant was improved to $1/1088$ by Yang and Zeng~\cite{YZ}, to $2/21$ by Gutin, Hao and Zhou~\cite{GHZ}, and to $1/5$ by Zeng~\cite{Zeng}.
	
	In this paper we count edges instead of vertices. For a graph $H$ we write $e(H)$ for the number of its edges. For a graph $G$, let
	\[
	\begin{aligned}
		f_o^*(G)&=\max\{e(G[S]):S\subseteq V(G),\ G[S]\text{ is odd}\},\\
		f_e^*(G)&=\max\{e(G[S]):S\subseteq V(G),\ G[S]\text{ is even}\}.
	\end{aligned}
	\]
	Every edge induces an odd subgraph, so $f_o^*(G)\ge1$ whenever $G$ has an edge. Scott~\cite[Section~6]{Scott01} asked whether there is a constant $c>0$ such that $f_o^*(G)\ge c\,e(G)$ for every graph $G$. A linear lower bound on $f_o(G)$ does not answer this question. Indeed, every vertex of an odd graph is incident with an edge, so the theorem of Zeng mentioned above gives $f_o^*(G)\ge n/10$, where $n$ is the number of non-isolated vertices of $G$; but $e(G)$ may be as large as $\binom n2$, and so this only yields $f_o^*(G)=\Omega(\sqrt{e(G)})$. {In a recent preprint, Liu~\cite{Liu} restated Scott's question and proved that $f_o^*(G)=\Omega(m^{4/7})$ for general graphs and $f_o^*(G)=\Omega(m^{2/3})$ for bipartite graphs, where $m=e(G)$.} Corollary \ref{cor:main}  answers the question of Scott in the affirmative.

	\begin{theorem}\label{thm:main}
		Let $G$ be a graph with $m$ edges and $n$ non-isolated vertices. Then
		\[
		m\le \YZ{1928f_o^*(G)+88n}.
		\]
	\end{theorem}
	
	As noted above, $n\le10f_o^*(G)$ by the theorem of Zeng, and so Theorem~\ref{thm:main} has the following consequence.
	
	\begin{corollary}\label{cor:main}
		Every graph $G$ with $m$ edges satisfies $f_o^*(G)\ge m/\YZ{2808}\ge m/\YZ{3000}$.
	\end{corollary}
	
	Any linear lower bound on \(f_o(G)\) can be used here in place of
Zeng’s theorem, possibly at the cost of a smaller coefficient of \(m\).
	
	Scott also noted in~\cite[Section~6]{Scott01} that the analogous question for $f_e^*(G)$ is open. As stated, its answer is negative: an induced subgraph of a forest is even only if it has no edges, so $f_e^*(G)=0$ for every forest $G$. A forest has fewer edges than vertices, and our second result shows that an even induced subgraph with many edges exists as soon as the number of edges is large compared with the number of vertices.
	
	\begin{theorem}\label{thm:even}
		Every graph $G$ with $n$ vertices and $m$ edges satisfies $f_e^*(G)\ge\YZ{(m-11n)/120}$.
	\end{theorem}
	
	We say that a set $S$ of vertices of a graph is \emph{odd} (\emph{even}, respectively) if it induces an odd (even, respectively) subgraph. Theorems~\ref{thm:main} and~\ref{thm:even} are both derived from the following result, in which the pairs are not required to be edges.
	
	\begin{theorem}\label{thm:pairs}
		Let $Q$ be a graph on $s$ vertices, and let $\cM$ be a set of unordered pairs of distinct vertices of $Q$. Then $Q$ has an even set that contains both vertices of at least \YZ{$(|\cM|-11s)/120$} pairs of $\cM$.
	\end{theorem}
	
	Theorem~\ref{thm:even} is the special case of Theorem~\ref{thm:pairs} in which $\cM$ is the edge set of the graph.
	
	The proofs are based on the following ideas. The even sets of a graph $Q$ are transformed in a simple way when a local complementation is applied to $Q$. Using this, we construct by induction on the number of vertices a random even set of $Q$ with the following property: for every non-empty set $U$ of vertices, the probability that the random even set contains an odd number of vertices of $U$ lies between $1/3$ and $2/3$.  The distribution also satisfies an entropy inequality. An entropy argument then shows that all pairs of vertices, apart from the edges of a \YZ{$11$}-degenerate graph (i.e., a graph in which every nonempty induced subgraph has a vertex of degree at most \YZ{11}) are contained in the random even set with probability at least \YZ{$1/120$}, which gives Theorem~\ref{thm:pairs}.
	
	To pass from even sets to odd sets we use a Schur complement over the field $\F$. For suitable disjoint vertex sets $P$ and $W$ of $G$, it gives a graph $Q$ on $W$ such that every even set of $Q$ becomes an odd set of $G$ after certain vertices of $P$ are added. The graph $Q$ is in general different from $G[W]$, so an edge of $G[W]$ need not be an edge of $Q$; this is the reason why Theorem~\ref{thm:pairs} is stated for arbitrary pairs. If $B$ is a maximal independent set of $G$ and $P$ is a subset of $B$, this approach bounds the number of edges of $G-B$. To count the edges between $B$ and its complement, we allow $P$ to contain adjacent vertices,
	and the Schur complement requires a certain matrix to be invertible. We guarantee this by means of private neighbours: $P$ consists of a subset of a set $C$ in which every vertex $c$ has a neighbour whose only neighbour in $C$ is $c$, together with some of these private neighbours. The edges incident with $C$ and with the private neighbours are counted separately.
	
	The paper is organised as follows. Section~\ref{sec:prelim} contains notation and preliminaries. In Section~\ref{sec:law} we construct the random even set and prove Theorems~\ref{thm:pairs} and~\ref{thm:even}. Section~\ref{sec:lift} contains the Schur complement and its consequences. Theorem~\ref{thm:main} and Corollary~\ref{cor:main} are proved in Section~\ref{sec:proof}, and Section~\ref{sec:remarks} contains concluding remarks.
	
	\section{Preliminaries}\label{sec:prelim}
	
	Let $G$ be a graph with vertex set $V(G)$ and edge set $E(G)$, and let $e(G)=|E(G)|$. For a vertex $v$, we write $N_G(v)$ for the set of neighbours of $v$, and $d_G(v)=|N_G(v)|$ for its degree. For $S\subseteq V(G)$, we write $G[S]$ for the subgraph of $G$ induced by $S$, $G-S$ for $G[V(G)\setminus S]$, and $d_S(v)=|N_G(v)\cap S|$ for the number of neighbours of $v$ in $S$; we write $G-v$ for $G-\{v\}$. For disjoint sets $S,T\subseteq V(G)$, we write $e(S,T)$ for the number of edges of $G$ with one end in $S$ and the other end in $T$. The subscript $G$ is omitted when the graph is clear from the context. In this notation, a set $S\subseteq V(G)$ is odd if $d_S(v)$ is odd for every $v\in S$, and even if $d_S(v)$ is even for every $v\in S$. Thus the empty set is both odd and even, and a vertex of an even set may have no neighbours in the set. For a positive integer $s$, let $[s]=\{1,\dots,s\}$.
	
	We use the following theorem of Zeng.
	
	\begin{theorem}[\cite{Zeng}]\label{thm:zeng}
		Every graph on $n$ vertices without isolated vertices has an odd set with at least $n/5$ vertices.
	\end{theorem}
	
	An odd set with $k$ vertices induces at least $k/2$ edges, as each of its vertices has a neighbour in the set. Deleting the isolated vertices of a graph $G$ changes neither $e(G)$ nor $f_o^*(G)$. Therefore Theorem~\ref{thm:zeng} gives $n\le10f_o^*(G)$ for every graph $G$ with $n$ non-isolated vertices.
	
	For a vertex $v$ of a graph $G$, the \emph{local complementation} of $G$ at $v$ is the graph $G*v$ obtained from $G$ by replacing the subgraph induced by $N_G(v)$ with its complement. That is, two distinct neighbours of $v$ are adjacent in $G*v$ if and only if they are not adjacent in $G$, and all other adjacencies are the same in $G*v$ and in $G$. In particular, $N_{G*v}(v)=N_G(v)$.
	
	All matrices and vectors are over the field $\F$ with two elements. Vectors are column vectors indexed by a finite set, usually a set of vertices. For a subset $S$ of the index set, $\one_S$ is the indicator vector of $S$, and $e_v=\one_{\{v\}}$. The \emph{support} $\supp x$ of a vector $x$ is the set of indices $v$ with $x_v=1$. We write $\A_G$ for the adjacency matrix of $G$ and, for $P,W\subseteq V(G)$, we write $\A_{P,W}$ for the submatrix of $\A_G$ with rows indexed by $P$ and columns indexed by $W$. If $x=\one_S$ for a set $S\subseteq V(G)$, then $(\A_Gx)_v$ is the parity of $d_S(v)$. Hence $S$ is even if and only if $(\A_Gx)_v=0$ for every $v\in S$, and $S$ is odd if and only if $(\A_Gx)_v=1$ for every $v\in S$. For a vector $\kappa$, we write $\Diag(\kappa)$ for the diagonal matrix with diagonal $\kappa$, and for a square matrix $K$, we write $\diag(K)$ for the vector of its diagonal entries. %If $K$ is a symmetric matrix and $z$ is a vector, then $z^\top Kz=\diag(K)^\top z$, because the terms $K_{ij}z_iz_j$ and $K_{ji}z_jz_i$ with $i\ne j$ cancel and $z_i^2=z_i$.
	
	A random subset $T$ of a finite set $U$ is \emph{uniform} if each element of $U$ belongs to $T$ independently with probability $1/2$.
	
	\begin{lemma}\label{lem:parity}
		Let $T$ be a uniform random subset of a finite set $U$, and let $N$ and $N'$ be non-empty subsets of $U$. Then $|N\cap T|$ is odd with probability $1/2$, and the probability that both $|N\cap T|$ and $|N'\cap T|$ are odd is at least $1/4$.
	\end{lemma}
	
	\begin{proof}
		Let $u\in N$. For each choice of $T\setminus\{u\}$, exactly one of the two choices for $u$ makes $|N\cap T|$ odd, which proves the first assertion and settles the case $N=N'$ of the second. If $N\ne N'$, we may assume that there is an element $u\in N\setminus N'$. The parity of $|N'\cap T|$ is determined by $T\setminus\{u\}$ and is odd with probability $1/2$, and for each choice of $T\setminus\{u\}$ the set $N\cap T$ has odd size with probability $1/2$. Hence both sizes are odd with probability $1/4$.
	\end{proof}

	%{\color{blue} \begin{proof}
			%Let $x=\one_T$. A linear map sends a uniform vector to a uniform vector on its image, since all non-empty fibres have the same size. Thus $\one_N^\top x$ is a fair bit. If $N=N'$, the parities coincide and are both odd with probability $1/2$. Otherwise $\one_N$ and $\one_{N'}$ are linearly independent over $\F$, so the pair of parities is uniform on $\F^2$ and equals $(1,1)$ with probability $1/4$.
			%\end{proof}}
			
			Logarithms are to base $2$. For a random variable $Z$ with finitely many values, $H(Z)=-\sum_z\Pr(Z=z)\log\Pr(Z=z)$ is its entropy, and for two such random variables, $H(Z\mid Z')=H(Z,Z')-H(Z')$ is the conditional entropy of $Z$ given $Z'$. If $H(Z\mid Z'=z')$ denotes the entropy of the conditional distribution of $Z$ given $Z'=z'$, then $H(Z\mid Z')=\sum_{z'}\Pr(Z'=z')H(Z\mid Z'=z')$. For $t\in[0,1]$, let $h(t)=-t\log t-(1-t)\log(1-t)$, with the convention $0\log 0=0$; this is the entropy of a random bit that equals $1$ with probability $t$, and $h$ is increasing on $[0,1/2]$. For probability distributions $\nu$ and $\mu$ on a finite set such that $\mu(z)>0$ whenever $\nu(z)>0$, the relative entropy is $D(\nu\,\|\,\mu)=\sum_z\nu(z)\log(\nu(z)/\mu(z))$, where the sum is over all $z$ with $\nu(z)>0$. We use the following standard facts; see~\cite{CT}.
			\begin{itemize}
				\item[(H1)] $H(Z_1,\dots,Z_k)=\sum_{i=1}^kH(Z_i\mid Z_1,\dots,Z_{i-1})$.
				\item[(H2)] $H(Z\mid Z',Z'')\le H(Z\mid Z')\le H(Z)$, and $H(Z\mid Z')=H(Z)$ if $Z$ and $Z'$ are independent. If $Z'$ is a function of $Z$, then $H(Z')\le H(Z)$, with equality if the function is injective.
				\item[(H3)] If $\Theta$ is a random variable with finitely many values, then $H(Z)-H(Z\mid\Theta)=\sum_\theta\Pr(\Theta=\theta)D(\nu_\theta\,\|\,\nu)$, where $\nu$ is the distribution of $Z$, $\nu_\theta$ is the conditional distribution of $Z$ given $\Theta=\theta$, and the sum is over all $\theta$ with $\Pr(\Theta=\theta)>0$.
				\item[(H4)] $D(\nu\,\|\,\mu)\ge0$. Moreover, if $\nu(C)=1$ for a set $C$, then $D(\nu\,\|\,\mu)\ge\log(1/\mu(C))$.
			\end{itemize}

			Finally, for vectors $g,x\in\F^V$ let $\chi_g(x)=(-1)^{g^\top x}$, where the exponent is read as an integer in $\{0,1\}$. If $X$ is a random vector in $\F^V$, then $\Pr(g^\top X=1)=(1-\EE\chi_g(X))/2$.
			
			\section{A random even set}\label{sec:law}
			
			Throughout this section, $Q$ is a graph with vertex set $V$ and $\A=\A_Q$. For $x\in\F^V$ and $v\in V$, the \emph{state} of $v$ in $x$ is the pair $(x_v,(\A x)_v)$. We say that $x$ is an \emph{even vector} of $Q$ if $x=\one_S$ for some even set $S$ of $Q.$ By the remarks in Section~\ref{sec:prelim}, the vector $x$ is even if and only if no vertex has state $(1,1)$ in $x$. If $x$ is an even vector, then we say that $v$ has \emph{type} $1$, $2$ or $0$ in $x$ if its state is $(1,0)$, $(0,1)$ or $(0,0)$, respectively. Thus $v$ has type $1$ if it belongs to the even set, type $2$ if it does not belong to the even set and has an odd number of neighbours in it, and type $0$ otherwise.
			
			The first lemma describes the effect of a local complementation on the states.
			
			\begin{lemma}\label{lem:lc}
				Let $j\in V$, let $r=\one_{N_Q(j)}$, and let $\A'$ be the adjacency matrix of $Q*j$. Define $\varphi\colon\F^V\to\F^V$ by $\varphi(x)=x+(1+r^\top x)e_j$. Then $\varphi$ is an affine bijection, and for every $x\in\F^V$ the vector $x'=\varphi(x)$ has the following properties.
				\begin{itemize}
					\item[\textup{(a)}] $x'_j=x_j+1+(\A x)_j$ and $(\A'x')_j=(\A x)_j$.
					\item[\textup{(b)}] If $i\in N_Q(j)$, then $x'_i=x_i$ and $(\A'x')_i=(\A x)_i+1+x_i$.
					\item[\textup{(c)}] If $i\notin N_Q(j)\cup\{j\}$, then $x'_i=x_i$ and $(\A'x')_i=(\A x)_i$.
				\end{itemize}
				In particular, $\varphi$ maps the set of even vectors of $Q$ onto the set of even vectors of $Q*j$.
			\end{lemma}
			
			\begin{proof}
				Note that $r^\top x=(\A x)_j$ and $r_j=0$. The map $x\mapsto x+(r^\top x)e_j$ is linear, and applying it twice gives the identity because $r^\top e_j=0$. Hence $\varphi$, which is this map followed by the translation by $e_j$, is an affine bijection. The vectors $x'$ and $x$ differ at most in the coordinate $j$, so $x'_i=x_i$ for every $i\ne j$, and $r^\top x'=r^\top x$. Since $x'_j=x_j+1+r^\top x$, this proves the statements about $x'$.
				
				The matrices $\A'$ and $\A$ differ exactly in the entries indexed by two distinct neighbours of $j$, and so $\A'=\A+rr^\top+\Diag(r)$. Moreover, $\A x'=\A x+(1+r^\top x)\A e_j$ and $(\A e_j)_i=r_i$. Therefore, for every $i\in V$,
				\[
				(\A'x')_i=(\A x)_i+r_i(1+r^\top x)+r_i\,r^\top x'+r_ix'_i=(\A x)_i+r_i(1+x'_i).
				\]
				For $i=j$ we have $r_j=0$, which gives (a). For $i\in N_Q(j)$ we have $r_i=1$ and $x'_i=x_i$, which gives (b), and for all other $i$ we have $r_i=0$, which gives (c).
				
				By (a), (b) and (c), the state of a vertex in $x'$ with respect to $Q*j$ is obtained from its state $(a,b)$ in $x$ with respect to $Q$ as follows: it is $(a+1+b,b)$ for the vertex $j$, it is $(a,b+1+a)$ for the neighbours of $j$, and it is $(a,b)$ for all other vertices. Each of these three maps is a permutation of $\F^2$ that fixes $(1,1)$. Hence $x$ is an even vector of $Q$ if and only if $x'$ is an even vector of $Q*j$. As $\varphi$ is a bijection, the last assertion follows.
			\end{proof}

			%{\color{blue} \begin{proof}
					%Put $b=r^\top x=(\A x)_j$. Changing $x_j$ leaves $b$ fixed, so adding $1+b$ to $x_j$ is an affine involution. The neighbourhood of $j$ is unchanged by local complementation, giving (a). For $i\in N_Q(j)$, changing $x_j$ contributes $1+b$ to the change in neighbour parity at $i$, and toggling its adjacencies to $N_Q(j)\setminus\{i\}$ contributes $b+x_i$. The total change is $1+x_i$, proving (b). All other states stay fixed, giving (c).
					
					%At $j$ the first bit changes only when the second is zero, and at a neighbour the second bit changes only when the first is zero. Thus a vertex has state $(1,1)$ after the transformation exactly when it had that state before. The involution therefore takes even vectors onto even vectors.
					%\end{proof}}
					
					For $y\in\F^{V\setminus\{v\}}$ we write $(y,0)$ for the vector in $\F^V$ that agrees with $y$ on $V\setminus\{v\}$ and has coordinate $0$ at $v$. The next lemma splits the even vectors of $Q$ into three families, according to the type of a fixed vertex $v$, and identifies the union of any two of the families with the set of even vectors of a graph on $V\setminus\{v\}$. The three graphs are obtained from $Q$ by deleting $v$ after at most two local complementations.
					
					\begin{lemma}\label{lem:three} 
						Let $v$ be a non-isolated vertex of $Q$, let $u$ be a neighbour of $v$, and let $r=\one_{N_Q(v)}$. Let
						\[
						\Lambda_0=\{x\in\F^V:x_v+r^\top x=1\},\quad
						\Lambda_1=\{x\in\F^V:x_v=0\},\quad
						\Lambda_2=\{x\in\F^V:r^\top x=0\},
						\]
						and let $Q_0=(Q*v)-v$, $Q_1=Q-v$ and $Q_2=((Q*u)*v)-v$.
						\begin{itemize}
							\item[\textup{(a)}] If $x$ is an even vector of $Q$ in which $v$ has type $t$, then $x\notin\Lambda_t$ and $x\in\Lambda_s$ for both $s\in\{0,1,2\}\setminus\{t\}$.
							\item[\textup{(b)}] For every $t\in\{0,1,2\}$ there is an affine bijection $\iota_t\colon\F^{V\setminus\{v\}}\to\Lambda_t$ that maps the set of even vectors of $Q_t$ onto the set of even vectors of $Q$ contained in $\Lambda_t$.
						\end{itemize}
					\end{lemma}
					\begin{proof}
						Write $(a,b)=(x_v,r^\top x)$. Among the equations $a+b=1$, $a=0$ and $b=0$, the state of type $t$ fails precisely the equation defining $\Lambda_t$. This proves (a).
						
						For (b), let $\psi_0,\psi_1,\psi_2$ be the affine maps induced by complementing at $v$, doing nothing, and complementing successively at $u$ and $v$, respectively. By Lemma~\ref{lem:lc} the resulting $v$-coordinates are $a+1+b$, $a$ and $b$. For the last calculation, the first complementation changes the state at $v$ to $(a,b+1+a)$, so the second gives $a+1+(b+1+a)=b$. Thus $\psi_t(x)_v=0$ exactly on $\Lambda_t$. Appending a zero at $v$ preserves evenness in both directions, as does each local complementation by Lemma~\ref{lem:lc}. Hence $\iota_t(y)=\psi_t^{-1}(y,0)$ has the required properties.
					\end{proof}

					Let $\gamma=\log(3/2)$. The following lemma provides the random even set described in the introduction.
					
					\begin{lemma}\label{lem:law}
						For every graph $Q$ with vertex set $V$ there is a random vector $X$ in $\F^V$ with the following properties.
						\begin{itemize}
							\item[\textup{(i)}] $X$ is always an even vector of $Q$.
							\item[\textup{(ii)}] $\Pr(X_v=1)=1/3$ for every $v\in V$.
							\item[\textup{(iii)}] $|\EE\chi_g(X)|\le1/3$ for every non-zero $g\in\F^V$.
							\item[\textup{(iv)}] $H(\Phi X)\ge\gamma\rank\Phi$ for every matrix $\Phi$ over $\F$ whose columns are indexed by $V$.
						\end{itemize}
					\end{lemma}
					
				\begin{proof}
					We induct on $|V|$. If $Q$ is edgeless, choose the coordinates of $X$ independently with $\Pr(X_v=1)=1/3$.  \YZ{Since every vector is even, (i) holds, and (ii) holds by construction. Moreover, by the independence of the coordinates of $X$, we have $\EE\chi_g(X)=(1/3)^{|\supp g|}$, which proves (iii).}

					\YZ{For a matrix $\Phi$ of rank $k$, let $I\subseteq V$ be a set of $k$ vertices such that the columns of $\Phi$ indexed by $I$ are linearly independent. For any fixed value of $(X_v)_{v\in V\setminus I}$, the vectors $\Phi X$ and $(X_v)_{v\in I}$ determine each other. Hence, by (H2) and the independence of the coordinates of $X$, we have
						\[
						\begin{aligned}
							H(\Phi X)&\ge H(\Phi X\mid (X_v)_{v\in V\setminus I})\\
							&=H((X_v)_{v\in I}\mid (X_v)_{v\in V\setminus I})\\
							&=H((X_v)_{v\in I})=kh(1/3)\ge\gamma k,
						\end{aligned}
						\]
						where the last inequality follows from $h(1/3)=\gamma+1/3$.}
					
					Otherwise choose an edge $uv$ and use the notation and maps of Lemma~\ref{lem:three}. For $t=0,1,2$,  \YZ{let $Y_t$ be a random vector with properties (i)--(iv) for $Q_t$} and put $X_t=\iota_t(Y_t)$. Thus $X_t$ is even and lies in $\Lambda_t$.  \YZ{Let $\nu_0=e_v+r$, $\nu_1=e_v$ and $\nu_2=r$.
						Then $\Lambda_t$ is the set of all $x\in\F^V$ for which
						$\nu_t^\top x$ takes a prescribed value.
						Since $r\ne0$ and $r_v=0$, the vectors $\nu_0,\nu_1,\nu_2$
						are non-zero and distinct.}  \YZ{Consequently, for every $g\notin\{0,\nu_t\}$, the linear form $g^\top x$ is non-constant on $\Lambda_t$. Thus, for every $g\notin\{0,\nu_t\}$, since $\iota_t$ is an affine bijection onto $\Lambda_t$, there exist a non-zero vector $g_t\in\F^{V\setminus\{v\}}$ and a constant $c_t\in\F$ such that $g^\top\iota_t(y)=g_t^\top y+c_t$ for every $y\in\F^{V\setminus\{v\}}$. By (iii) for $Y_t$, we have
						\begin{equation}\label{eq:indh}
								|\EE\chi_g(X_t)|=|\EE\chi_{g_t}(Y_t)|\le1/3.
						\end{equation}
						Consequently, each value of $g^\top X_t$ has probability between $1/3$ and $2/3$.}
					
					Let $\pi_t(s)$ be the probability that $v$ has type $s$ in $X_t$. Lemma~\ref{lem:three}(a) gives $\pi_t(t)=0$. For $s\ne t$, type $s$ is the event $X_t\notin\Lambda_s$, so $\pi_t(s)\in[1/3,2/3]$ by the preceding bound with $g=\nu_s$.

					\YZ{We now show that there exist non-negative real numbers $\lambda_0,\lambda_1,\lambda_2$ with sum $1$ such that $\sum_{t=0}^2\lambda_t\pi_t(s)=1/3$ for every $s\in\{0,1,2\}$. Since $\pi_0(0)=0$ and $\pi_1(0),\pi_2(0)\ge1/3$, for each $j\in\{1,2\}$ there is a convex combination $p_j$ of $\pi_0$ and $\pi_j$ such that $p_j(0)=1/3$. Since $\pi_0(2),\pi_1(2)\ge1/3$, we have $p_1(2)\ge1/3$. Similarly, since $\pi_0(1),\pi_2(1)\ge1/3$, we have
					$p_2(1)\ge1/3$ and hence $p_2(2)\le1/3$. Thus there is a convex combination $p$ of $p_1$ and $p_2$ such that $p(2)=1/3$. As $p(0)=1/3$, we also have $p(1)=1/3$. This proves the claim.}
					
					Choose $\Theta=t$ with probability $\lambda_t$ and, conditional on this choice, draw $X$ with the law of $X_t$.  \YZ{Since each $X_t$ is an even vector of $Q$ by
						Lemma~\ref{lem:three}(b), $X$ is always an even vector
						of $Q$, which proves (i).
						By the choice of $\lambda_0,\lambda_1,\lambda_2$,
						the vertex $v$ has each type in $X$ with probability $1/3$.} In particular\YZ{, by Lemma~\ref{lem:three}(a),} $\Pr(X\in\Lambda_t)=2/3$ for every $t$.
					
					\YZ{We now prove (iii). Since $\Pr(X\in\Lambda_t)=2/3$
						for every $t\in\{0,1,2\}$, by the definitions of
						$\Lambda_0,\Lambda_1,\Lambda_2$, we have
						$\Pr(\nu_0^\top X=1)=2/3$ and
						$\Pr(\nu_1^\top X=0)=\Pr(\nu_2^\top X=0)=2/3$.
						Hence
						\[
						\EE\chi_{\nu_0}(X)=1-2\cdot\frac23=-\frac13
						\quad\text{and}\quad
						\EE\chi_{\nu_1}(X)=\EE\chi_{\nu_2}(X)
						=1-2\cdot\frac13=\frac13.
						\]
						Thus $|\EE\chi_g(X)|=1/3$ for every
						$g\in\{\nu_0,\nu_1,\nu_2\}$.}
					
					\YZ{For every non-zero $g\notin\{\nu_0,\nu_1,\nu_2\}$, by \eqref{eq:indh}, we have
						\[
						|\EE\chi_g(X)|
						=\left|\sum_t\lambda_t\EE\chi_g(X_t)\right|
						\le\sum_t\lambda_t/3=1/3.
						\]
						This proves (iii).}
					
					To see (ii), take any non-isolated $w$ and put $f=\one_{N_Q(w)}$. Since $X$ is even, the three types at $w$ are distinguished by the disjoint events
					\[
					X_w=1,\qquad f^\top X=1,\qquad (e_w+f)^\top X=0.
					\]
				
				\YZ{Since $w$ is not isolated, we have $f\ne0$ and $f_w=0$, so the vectors $e_w$, $f$ and $e_w+f$ are all non-zero. By (iii), for every non-zero vector $g$, we have $|1-2\Pr(g^\top X=1)|\le1/3$, and hence each of the two possible values of $g^\top X$ has probability at least $1/3$. Thus each of the three events above has probability at least $1/3$. Since these events are pairwise disjoint and their probabilities sum to $1$, each has probability exactly $1/3$. In particular, $\Pr(X_w=1)=1/3$.}
				
				An isolated $w$ stays isolated in each $Q_t$. The maps $\iota_t$ leave its coordinate unchanged\YZ{, since the local complementations used in the construction of these maps are performed only at $u$ and $v$}\YZ{. Thus, by the induction hypothesis, $\Pr((X_t)_w=1)=\Pr((Y_t)_w=1)=1/3$ for every $t$, and hence $\Pr(X_w=1)=\sum_t\lambda_t\Pr((X_t)_w=1)=1/3$. This proves (ii).}
					
					For (iv), \YZ{let} $Z=\Phi X$, $k=\rank\Phi$, and  $a_t$ equal $1$ if $\YZ{\nu_t^{\top}}$ belongs to the row space of $\Phi$ and $0$ otherwise.

					\YZ{Since $\iota_t$ is affine, we can write $\iota_t(y)=L_ty+b_t$, where $L_t$ is a matrix and $b_t$ is a vector. Since $\iota_t$ is an affine bijection onto $\Lambda_t$, the linear map $L_t$ is an isomorphism from $\F^{V\setminus\{v\}}$ onto $\nu_t^\perp=\{x\in\F^V:\nu_t^\top x=0\}$. Consider the linear map $g^\top\mapsto g^\top L_t$ from the row space of $\Phi$ to the row space of $\Phi L_t$. Since the image of $L_t$ is $\nu_t^\perp$, we have $g^\top L_t=0$ if and only if $g\in\{0,\nu_t\}$. Hence the kernel of this map has dimension $a_t$, and therefore $\rank(\Phi L_t)=k-a_t$ by the rank-nullity theorem.
						
						Conditional on $\Theta=t$, the random vector $Z$ has the same distribution as $\Phi\iota_t(Y_t)=\Phi L_tY_t+\Phi b_t$. Since translation does not change entropy, by (iv) for $Y_t$ we have, whenever $\lambda_t>0$,
						\[
						H(Z\mid\Theta=t)=H(\Phi L_tY_t)
						\ge\gamma\rank(\Phi L_t)=\gamma(k-a_t).
						\]}

					\YZ{Let $\mu$ be the distribution of $Z$, and for every $t$ with $\lambda_t>0$, let $\mu_t$ be the conditional distribution of $Z$ given $\Theta=t$. We claim that $D(\mu_t\,\|\,\mu)\ge\gamma a_t$. If $a_t=0$, this follows from the non-negativity of relative entropy. Suppose that $a_t=1$. Since $\nu_t^\top$ belongs to the row space of $\Phi$, there is a vector $z$ such that $\nu_t=\Phi^\top z$. Hence $\nu_t^\top x=z^\top\Phi x$ for every $x\in\F^V$, and therefore $\Phi x$ determines whether $x$ belongs to $\Lambda_t$.
						
						Let $C_t=\{\Phi x:x\in\Lambda_t\}$. By the above observation, $\Phi x\in C_t$ if and only if $x\in\Lambda_t$. Consequently, $\mu(C_t)=\Pr(X\in\Lambda_t)=2/3$. On the other hand, conditional on $\Theta=t$, we have $X\in\Lambda_t$ with probability $1$, so $\mu_t(C_t)=1$. By (H4), we obtain $D(\mu_t\,\|\,\mu)\ge\log(1/\mu(C_t))=\gamma$, proving the claim.}
					
					Now (H3) gives
					
					\YZ{
						\[
						\begin{aligned}
							H(Z)&=\sum_{t:\lambda_t>0}\lambda_t
							\bigl(H(Z\mid\Theta=t)+D(\mu_t\,\|\,\mu)\bigr)\\
							&\ge\sum_{t:\lambda_t>0}\lambda_t
							\bigl(\gamma(k-a_t)+\gamma a_t\bigr)
							=\gamma k.
						\end{aligned}
						\]
						This proves (iv) and completes the proof of the lemma.}
				\end{proof}
					
					Properties (iii) and (iv) imply that most pairs of coordinates of $X$ are equal to $1$ simultaneously with probability bounded away from zero. For the application in Section~\ref{sec:proof} we need this not only for the coordinates of $X$ but for linearly independent linear forms in these coordinates.
					
					\begin{theorem}\label{thm:forms}
						Let $Q$ be a graph with vertex set $V$, and let $X$ be a random vector in $\F^V$ with properties (iii) and (iv) of Lemma~\ref{lem:law}. Let $g_1,\dots,g_s\in\F^V$ be linearly independent, let $c_1,\dots,c_s\in\F$, and let $Z_i=g_i^\top X+c_i$ for $i\in[s]$. Let $\cM$ be a set of unordered pairs of distinct elements of $[s]$. Then the expected number of pairs $\{i,j\}\in\cM$ with $Z_i=Z_j=1$ is at least \YZ{$(|\cM|-11s)/120$}.
					\end{theorem}
					
					\begin{proof}
						\YZ{Let} $Z_S=(Z_i)_{i\in S}$ \YZ{and $D$ be a non-empty subset of $[s]$. Then, we have
						\[\sum_{j\in D}Z_j=g_D^\top X+c_D,\]
						 where $g_D=\sum_{j\in D}g_j$ and $c_D=\sum_{j\in D}c_j$. Since the vectors $g_1,\dots,g_s$ are linearly independent, we have $g_D\ne0$. By Lemma~\ref{lem:law}(iii),
							\begin{equation}\label{eq:fourier-bound}
									\Bigl|\EE(-1)^{\sum_{j\in D}Z_j}\Bigr|=|\EE\chi_{g_D}(X)|\le1/3.
							\end{equation}
							In particular, when $D=\{i\}$, it gives $\Pr(Z_i=1)\ge1/3$ for every $i\in[s]$. Moreover, let $\Phi_D$ be the matrix whose rows are $g_i^\top$, $i\in D$. Since adding a fixed vector is a bijection, by (H2) and (iv) of Lemma~\ref{lem:law}, we have
							\begin{equation}\label{eq:entropylb}
									H(Z_D)=H(\Phi_D X)\ge\gamma\rank\Phi_D=\gamma|D|.
							\end{equation}
						}
						
						\YZ{Construct a graph $F$ with vertex set $[s]$ by adding the edge $ij$} when $\Pr(Z_i=Z_j=1)<\YZ{1/120}$.

						\YZ{For $1\le k\le5$, let $B$ be a set of $k$ neighbours of a vertex $i$ in $F$, and let $E$ be the event that $Z_j=0$ for every $j\in B$. Since the indicator of $E$ is $\prod_{j\in B}\frac12(1+(-1)^{Z_j})$, by \eqref{eq:fourier-bound}, we have
							\[
							\Pr(E)=\frac1{2^k}\sum_{D\subseteq B}\EE(-1)^{\sum_{j\in D}Z_j}
							\le\frac{1+(2^k-1)/3}{2^k}
							=\frac13+\frac{2}{3\cdot2^k}.
							\]
							Moreover, since every $j\in B$ is adjacent to $i$ in $F$, we have
							$\Pr(Z_i=1,E^c)\le\sum_{j\in B}\Pr(Z_i=Z_j=1)\le k/120$.
							
							Let $\widehat Z_i$ be the indicator of $E$. As $\widehat Z_i\ne Z_i$ if and only if either $E$ occurs and $Z_i=0$ or $E^c$ occurs and $Z_i=1$, we have
							\[
								p:=\Pr(\widehat Z_i\ne Z_i)
								=\Pr(E)-\Pr(Z_i=1)+2\Pr(Z_i=1,E^c)
								\le\frac{2}{3\cdot2^k}+\frac{k}{60}.
							\]
							Since $\widehat Z_i$ is determined by $Z_B$, the bits $Z_i$ and $Z_i+\widehat Z_i$ determine each other given $Z_B$. Thus, by (H2) and the monotonicity of $h$ on $[0,1/2]$, we have
							\[
							H(Z_i\mid Z_B)=H(Z_i+\widehat Z_i\mid Z_B)
							\le H(Z_i+\widehat Z_i)=h(p)
							\le h\left(\frac{2}{3\cdot2^k}+\frac{k}{60}\right),
							\]
						where the last inequality follows from the fact that $\frac{2}{3\cdot2^k}+\frac{k}{60}<1/2$ for $1\le k\le5$.
							
					Let $\alpha_0=1$ and, for $1\le k\le5$, let
					$\alpha_k=h\left(\frac{2}{3\cdot2^k}+\frac{k}{60}\right)$.
					A computer-assisted computation gives
					\[
					\alpha_1\le0.935,\quad
					\alpha_2\le0.722,\quad
					\alpha_3\le0.567,\quad
					\alpha_4\le0.495,\quad
					\alpha_5\le0.483.
					\]
						By the above inequality and the fact that conditioning cannot
						increase entropy, for any set $B$ of $k$ neighbours of $i$ in $F$,
						we have 
					\[
							H(Z_i\mid Z_B)\le\alpha_{\min\{k,5\}}.
		\]
						}

						\YZ{Suppose that $F[S]$ has minimum degree at least $12$ for some $S\subseteq [s]$. Let $\sigma$ be a uniform random linear order of $S$. For every $i\in S$, let $P_i$ be the set of vertices of $S$ preceding $i$ in $\sigma$, let $d_i$ be the degree of $i$ in $F[S]$, and let $k_i$ be the number of neighbours of $i$ in $P_i$ in $F[S]$. Then $k_i$ is uniformly distributed on $\{0,\dots,d_i\}$, since the position of $i$ relative to its $d_i$ neighbours is uniform. By (H2) and the above bound, $H(Z_i\mid Z_{P_i})\le\alpha_{\min\{k_i,5\}}$ for every fixed order $\sigma$. Thus, applying (H1) to each order $\sigma$ and averaging over $\sigma$, we obtain
			\[
			\begin{aligned}
				H(Z_S)
				&=\sum_{i\in S}\EE_\sigma[H(Z_i\mid Z_{P_i})]\\
				&=\sum_{i\in S}\sum_{t=0}^{d_i}
				\Pr(k_i=t)\EE_\sigma[H(Z_i\mid Z_{P_i})\mid k_i=t]\\
				&\le\sum_{i\in S}\frac{1}{d_i+1}
				\sum_{t=0}^{d_i}\alpha_{\min\{t,5\}}\\
				&=\sum_{i\in S}\frac{1}{d_i+1}
				\left(\sum_{t=0}^4\alpha_t+(d_i-4)\alpha_5\right)\\
				&\le\frac{\sum_{t=0}^4\alpha_t+8\alpha_5}{13}|S|
				<0.584|S|<\gamma|S|,
			\end{aligned}
			\]
							where the second inequality follows from $d_i\ge12$ and $\alpha_k\ge\alpha_5$ for $0\le k\le4$. This contradicts the lower bound $H(Z_S)\ge\gamma|S|$. Hence $F$ is $11$-degenerate and has at most $11s$ edges. Each pair of $\cM$ outside $E(F)$ contributes at least $1/120$ to the expectation, and so the expected number of pairs $\{i,j\}\in\cM$ with $Z_i=Z_j=1$ is at least $(|\cM|-11s)/120$.}
					\end{proof}

					\YZ{\begin{proof}[Proof of Theorem~\ref{thm:pairs}]
						Let $V$ be the vertex set of $Q$, and let $X$ be a random
						vector as in Lemma~\ref{lem:law}. We apply
						Theorem~\ref{thm:forms} with $g_v=e_v$ and $c_v=0$
						for every $v\in V$. Since the vectors $e_v$, $v\in V$,
						are linearly independent and $g_v^\top X+c_v=X_v$,
						the expected number of pairs $\{v,w\}\in\cM$ with
						$X_v=X_w=1$ is at least $(|\cM|-11s)/120$.
						Hence there is an even vector $x$ of $Q$ such that
						the even set $\supp x$ contains both vertices of at least
						$(|\cM|-11s)/120$ pairs of $\cM$.
					\end{proof}}

					%{\color{blue} \begin{proof}[Proof of Theorem~\ref{thm:pairs}]
							%Take $X$ from Lemma~\ref{lem:law} and apply Theorem~\ref{thm:forms} to the coordinate forms $g_v=e_v$, with $c_v=0$. Since $\supp X$ is always even, some even set contains both ends of at least the expected number of pairs of $\cM$.
							%\end{proof}}
							
							\begin{proof}[Proof of Theorem~\ref{thm:even}]
								By Theorem~\ref{thm:pairs} applied to $Q=G$ and $\cM=E(G)$, there is an even set $S$ of $G$ \YZ{such that $G[S]$} contains  at least \YZ{$(m-11n)/120$ edges}. 
							\end{proof}
							
							\section{From even sets to odd sets}\label{sec:lift}
							
							In this section $G$ is a graph and $\A=\A_G$. The first lemma uses a Schur complement to relate even sets to odd sets.

							\begin{lemma}\label{lem:schur}
								Let $P$ and $W$ be disjoint sets of vertices of $G$, and let $M=\A_{P,W}$. Let $K$ be an invertible symmetric matrix with rows and columns indexed by $P$ whose off-diagonal entries are equal to those of $\A_{P,P}$, and let $\kappa=\diag(K)$. Suppose that $M^\top K^{-1}\one_P=\one_W$. Then there is a graph $Q$ with vertex set $W$ such that, for every even set $Y$ of $Q$, the set
								\[
								Y\cup\supp\bigl(K^{-1}(\one_P+\kappa+M\one_Y)\bigr)
								\]
								is an odd set of $G$.
							\end{lemma}
							\begin{proof}
								Let $\widehat\A=\A_{W,W}+M^\top K^{-1}M$ and $\delta=\diag(\widehat\A)$. The matrix $\widehat\A+\Diag(\delta)$ is symmetric and has zero diagonal, so it is the adjacency matrix of a graph $Q$ with vertex set $W$.
								
								%We first show that $\delta=M^\top K^{-1}\kappa$. Let $p\in P$ and $z=K^{-1}e_p$. By the identity $z^\top Kz=\diag(K)^\top z$ from Section~\ref{sec:prelim}, we have $e_p^\top K^{-1}e_p=z^\top Kz=\kappa^\top K^{-1}e_p$, which is the coordinate of $K^{-1}\kappa$ at $p$ because $K^{-1}$ is symmetric. Hence $\diag(K^{-1})=K^{-1}\kappa$. Now let $w\in W$, and let $c$ be the column of $M$ indexed by $w$. Applying the same identity to the symmetric matrix $K^{-1}$ and the vector $c$, we obtain $c^\top K^{-1}c=\diag(K^{-1})^\top c=\kappa^\top K^{-1}c$. The left-hand side is the diagonal entry of $M^\top K^{-1}M$ at $w$, and the right-hand side is the coordinate of $M^\top K^{-1}\kappa$ at $w$. Since $\A_{W,W}$ has zero diagonal, this shows that $\delta=M^\top K^{-1}\kappa$.
								
								We first show that $\delta=M^\top K^{-1}\kappa$. Since $K$ is symmetric and we are working over $\F$, for every $z\in\F^P$ we have
								\begin{equation}\label{eq:KinF2}
									z^\top Kz=\sum_i K_{ii}z_i^2+2\sum_{i<j}K_{ij}z_iz_j
									=\sum_i K_{ii}z_i=\kappa^\top z.
								\end{equation}
								Thus, for any column $c$ of $M$, letting $z=K^{-1}c$, we have
								\[
								c^\top K^{-1}c=(Kz)^\top z=z^\top Kz
								=\kappa^\top z=\kappa^\top K^{-1}c.
								\]
								Since $\A_{W,W}$ has zero diagonal and $K^{-1}$ is symmetric, this implies that $\delta=M^\top K^{-1}\kappa$.
								
								Let $Y$ be an even set of $Q$, let $y=\one_Y$ and $x=K^{-1}(\one_P+\kappa+My)$, and let $S=Y\cup\supp x$. We have to show that every vertex of $S$ has an odd number of neighbours in $S$. For the vertices of $P$, these parities are the coordinates of the vector $\A_{P,P}x+My$. Since $\A_{P,P}=K+\Diag(\kappa)$, we have
								\[
								\A_{P,P}x+My=Kx+\Diag(\kappa)x+My=\one_P+\kappa+\Diag(\kappa)x.
								\]
								At a vertex $p$ with $x_p=1$, the coordinate of this vector is $1+\kappa_p+\kappa_p=1$. For the vertices of $W$, the parities are the coordinates of the vector $\A_{W,W}y+M^\top x$. By the assumption $M^\top K^{-1}\one_P=\one_W$ and 
								the fact that $\delta=M^\top K^{-1}\kappa$,
								\[
								M^\top x=M^\top K^{-1}\one_P+M^\top K^{-1}\kappa+M^\top K^{-1}My=\one_W+\delta+M^\top K^{-1}My,
								\]
								and therefore
								\[
								\A_{W,W}y+M^\top x=\one_W+\delta+\widehat\A y=\one_W+\delta+\A_Qy+\Diag(\delta)y.
								\]
								At a vertex $w\in Y$ we have $(\A_Qy)_w=0$, because $Y$ is an even set of $Q$, and the coordinate of $\Diag(\delta)y$ is $\delta_w$. So the coordinate of the last vector at $w$ is $1+\delta_w+\delta_w=1$. Thus $S$ is an odd set of $G$.
							\end{proof}
							
							When $P$ is independent, we may take $K$ to be the identity matrix in Lemma~\ref{lem:schur}.
							
							\begin{lemma}\label{lem:indep}
								Let $P$ be an independent set of $G$, and let $W\subseteq V(G)\setminus P$ be a set such that $d_P(w)$ is odd for every $w\in W$. Then there is a graph $Q$ with vertex set $W$ such that, for every even set $Y$ of $Q$, the set $Y\cup\{p\in P:d_Y(p)\text{ is odd}\}$ is an odd set of $G$.
							\end{lemma}
							
							\begin{proof}
								Let $K$ be the identity matrix indexed by $P$. Its off-diagonal entries are equal to those of $\A_{P,P}$ because $P$ is independent, and $\kappa=\one_P$. Since $M=\A_{P,W}$, the coordinate of $M^\top K^{-1}\one_P=M^\top\one_P$ at $w\in W$ is the parity of $d_P(w)$, which is $1$. Hence Lemma~\ref{lem:schur} applies. For an even set $Y$ of the resulting graph $Q$, we have $K^{-1}(\one_P+\kappa+M\one_Y)=M\one_Y$, and therefore the support of this vector is the set of all $p\in P$ for which $d_Y(p)$ is odd.
							\end{proof}

							The next lemma allows us to use Theorem~\ref{thm:pairs} to bound the number of edges in certain induced subgraphs of $G$.
							
							\begin{lemma}\label{lem:count}
								Let $W\subseteq V(G)$, and suppose that there is a graph $Q$ with vertex set $W$ such that every even set of $Q$ is contained in an odd set of $G$. Then $e(G[W])\le\YZ{120}f_o^*(G)+\YZ{11}|W|$.
							\end{lemma}
							\begin{proof}
								By Theorem~\ref{thm:pairs} applied to $Q$ and $\cM=E(G[W])$, there is an even set $Y$ of $Q$ that contains both ends of at least $(e(G[W])-\YZ{11|W|)/120}$ edges of $G[W]$. Let $S$ be an odd set of $G$ with $Y\subseteq S$. All these edges belong to $G[S]$, and so their number is at most $f_o^*(G)$.
							\end{proof}

							In the following lemma, the set $P$ of Lemma~\ref{lem:schur} need not be independent. We enlarge a set $T$ of vertices by its private neighbours and choose the diagonal entries of the resulting matrix to ensure its invertibility.

							\begin{lemma}\label{lem:witness}
								Let $C\subseteq D\subseteq V(G)$, and suppose that for every $c\in C$ there is a vertex $w_c\in D$ such that $N_G(w_c)\cap C=\{c\}$. Let $T\subseteq C$, and let $W_T$ be the set of all vertices $v\in V(G)\setminus D$ for which $d_T(v)$ is odd. Then there is a graph $Q$ with vertex set $W_T$ such that every even set $Y$ of $Q$ is contained in an odd set $S$ of $G$ with $S\setminus Y\subseteq D$.
							\end{lemma}
							\begin{proof}
								If $T=\emptyset$, then $W_T=\emptyset$ and there is nothing to prove. Let $R=\{w_t:t\in T\}\setminus T$, $P=T\cup R$, and let $\eta\in\F^P$ be the indicator vector of $T$. Among all symmetric matrices indexed by $P$ whose off-diagonal entries are equal to those of $\A_{P,P}$ and which satisfy $K\eta=\one_P$, let $K$ be one of maximum rank. Such a matrix exists since every vertex of $R$ has exactly one neighbour in $T$, and we can choose the diagonal entries indexed by $T$ so that $(K\eta)_t=1$ for every $t\in T$.
								
								We claim that $K$ is invertible. Suppose otherwise, and let $z$ be a non-zero vector in $\ker K$. If $z_r=0$ for every $r\in R$, then for every $t\in T$ we have $(Kz)_{w_t}=z_t$. Indeed, $w_t$ has exactly one neighbour in $T$, namely $t$, and if $w_t\in T$, then $K_{w_t,w_t}=0$ since $(K\eta)_{w_t}=1$. As $Kz=0$, we obtain $z_t=0$ for every $t\in T$, contradicting that $z$ is non-zero.
								
								Hence there is a vertex $r\in R$ such that $z_r=1$. Let $K'=K+e_re_r^\top$. Since $\eta_r=0$, we have $K'\eta=K\eta=\one_P$. For every $u\in\ker K'$, by the symmetry of $K$ we have
								\[
								0=z^\top K'u=z^\top Ku+z_ru_r=u_r.
								\]
								Hence $Ku=K'u+u_re_r=0$, and so $\ker K'\subseteq\ker K$.
								Moreover, $z\notin\ker K'$ since $K'z=e_r\ne 0$.
								Thus $\operatorname{rank}K'>\operatorname{rank}K$,
								contradicting the choice of $K$. This proves that $K$ is invertible.
								
								Let $M=\A_{P,W_T}$. Since $K^{-1}\one_P=\eta$ and every vertex of $W_T$ has an odd number of neighbours in $T$, we have $M^\top K^{-1}\one_P=M^\top\eta=\one_{W_T}$. Hence Lemma~\ref{lem:schur} applies and gives a graph $Q$ with vertex set $W_T$ such that every even set $Y$ of $Q$ is contained in an odd set $S$ of $G$ with $S\setminus Y\subseteq P\subseteq D$.
							\end{proof}

							\begin{lemma}\label{lem:outside}
								Let $C$, $D$ and the vertices $w_c$ be as in Lemma~\ref{lem:witness}, and let $V_0\subseteq V(G)\setminus D$ be a set such that every vertex of $V_0$ has a neighbour in $C$. Then $e(G[V_0])\le\YZ{480f_o^*(G)+22|V_0|}$.
							\end{lemma}
							
							\begin{proof}
								Let $T$ be a uniform random subset of $C$, let $W_T$ be as in Lemma~\ref{lem:witness}, and let $U_T=V_0\cap W_T$. For an edge $vv'$ of $G[V_0]$, the sets $N_G(v)\cap C$ and $N_G(v')\cap C$ are non-empty, so by Lemma~\ref{lem:parity} both $v$ and $v'$ belong to $U_T$ with probability at least $1/4$. Moreover, every vertex of $V_0$ belongs to $U_T$ with probability $1/2$. Hence $\mathbb E[e(G[U_T])]\ge e(G[V_0])/4$ and $\mathbb E[|U_T|]=|V_0|/2$.
								
								By Lemma~\ref{lem:witness}, for every $T\subseteq C$ there is a graph $J$ with vertex set $W_T$ such that every even set of $J$ is contained in an odd set of $G$. Since every even set of $J[U_T]$ is also an even set of $J$, Lemma~\ref{lem:count} gives $e(G[U_T])\le\YZ{120f_o^*(G)+11|U_T|}$. Taking expectations, we obtain
								\[
								e(G[V_0])/4\le\YZ{120f_o^*(G)+11|V_0|/2},
								\]
								which completes the proof.
							\end{proof}

							\section{Proof of Theorem~\ref{thm:main}}\label{sec:proof}
							
							In the proof of Theorem~\ref{thm:main}, Lemma~\ref{lem:outside} will be applied to the set of all vertices outside a certain set $D$, and so it does not count the edges that have an end in $D$. These edges are counted by the following lemma. 
							
							\begin{lemma}\label{lem:matched}
								Let $V(G)=A\cup B$ be a partition of the vertex set of $G$ such that $B$ is independent. Let $C\subseteq A$, and for every $c\in C$ let $\ell_c\in B$ be a vertex such that $N_G(\ell_c)\cap C=\{c\}$. Let $L=\{\ell_c:c\in C\}$. Then the number of edges of $G$ between $A$ and $B$ that have an end in $C\cup L$ is at most $\YZ{484f_o^*(G)+22|A|}+|C|$.
							\end{lemma}

							\begin{proof}
								The edges between $C$ and $L$ form a matching of size $|C|$. We bound the other two contributions.
								
								Choose a uniform $U\subseteq C$ and put $B_U=\{b\in B\setminus L:d_U(b)\text{ is odd}\}$. Every vertex of $B_U$ has odd degree in $G[U\cup B_U]$ because $B$ is independent. For each $c\in U$ of even degree there, add $\ell_c$. Each added vertex has just one neighbour in the resulting set, namely $c$, so the resulting set is odd. For an edge $cb$ between $C$ and $B\setminus L$, Lemma~\ref{lem:parity} applied to $\{c\}$ and $N(b)\cap C$ gives $\Pr(c\in U,\,b\in B_U)\ge1/4$. Taking expectations gives $e(C,B\setminus L)\le4f_o^*(G)$.
								
								Next choose a uniform $T\subseteq L$ and put $W=\{a\in A:d_T(a)\text{ is odd}\}$. The same argument, now with $\{\ell\}$ and $N(a)\cap L$ for each edge $\ell a$, gives
								\begin{equation}\label{eq:l51}
									\frac14 e(L,A\setminus C)\le\EE\,e(T,W\setminus C).
								\end{equation}

									Fix $T$ and apply Lemma~\ref{lem:indep} with $P=T$ and $W$ as defined above. Since $T$ is independent and $d_T(a)$ is odd for every $a\in W$, there is a graph $Q$ with vertex set $W$ such that, for every even set $Y$ of $Q$, the set
									\[
									S_Y=Y\cup\{\ell\in T:d_Y(\ell)\text{ is odd}\}
									\]
									is an odd set of $G$.
									
									Let $X\in\F^W$ be a random vector obtained by applying Lemma~\ref{lem:law} to $Q$, and let $Y=\supp X$ (recall that by Lemma~\ref{lem:law}, $Y$ is always an even set of $Q$). For each $a\in W\setminus C$, let $g_a=e_a$, and for each $\ell\in T$, let $g_\ell=\one_{N_G(\ell)\cap W}$. Note that a vertex $v\in T\cup(W\setminus C)$ belongs to $S_Y$ if and only if $g_v^\top X=1$.
									
									We claim that the vectors $g_v$, $v\in T\cup(W\setminus C)$, are linearly independent. Indeed, for each $\ell_c\in T$, we have $c\in W$ since $d_T(c)=1$. Moreover, the coordinate of $g_{\ell_c}$ at $c$ is $1$, while the coordinate of every other vector in the family at $c$ is $0$. Hence, in any linear combination of these vectors that equals zero, the coefficient of every $g_\ell$, $\ell\in T$, must be zero. Since the remaining vectors $g_a$, $a\in W\setminus C$, are distinct unit vectors, their coefficients must also be zero. This proves the claim.

								Let $\cM$ be the set of all pairs $\{\ell,a\}$ with $\ell\in T$,
									$a\in W\setminus C$ and $\ell a\in E(G)$.
									Since $|T|=|W\cap C|$, the number of vectors $g_v$ is
									$|T|+|W\setminus C|=|W|$.
									Applying Theorem~\ref{thm:forms} to $X$, the linearly independent
									vectors $g_v$, $v\in T\cup(W\setminus C)$, and the set $\cM$
									(with all constant terms equal to zero), we obtain that the
									expected number of pairs $\{\ell,a\}\in\cM$ with
									$g_\ell^\top X=g_a^\top X=1$ is at least \YZ{$(|\cM|-11|W|)/120$}.
									Since every such pair is an edge of $G[S_Y]$ and $S_Y$ is
									an odd set of $G$, we have
									\[
									f_o^*(G)\ge \YZ{(|\cM|-11|W|)/120
									= (e(T,W\setminus C)-11|W|)/120,}
									\]
									and therefore \YZ{$e(T,W\setminus C)\le120	f_o^*(G)+11|W|$}.
									Moreover, by Lemma~\ref{lem:parity}, each vertex of $A$
									belongs to $W$ with probability at most $1/2$,
									and so $\EE|W|\le|A|/2$.
									Taking expectations, we obtain from \eqref{eq:l51} that
									\[
									e(L,A\setminus C)/4
									\le\YZ{120f_o^*(G)+11|A|/2.}
									\]
									Adding the three bounds completes the proof.
							\end{proof}

							\begin{proof}[Proof of Theorem~\ref{thm:main}]
								We may assume that $G$ has no isolated vertices, since deleting isolated vertices does not change any of the quantities in the statement. Thus $n=|V(G)|$. Let $B$ be a maximal independent set and $A=V(G)\setminus B$. Thus, every vertex has a neighbour in the other part. For a uniform $T\subseteq B$, put $W=\{a\in A:d_T(a)\text{ is odd}\}$. By Lemma~\ref{lem:parity}, both ends of each edge of $G[A]$ lie in $W$ with probability at least $1/4$. Moreover, every vertex of $A$ belongs to $W$ with probability $1/2$. Hence $\EE[e(G[W])]\ge e(G[A])/4$ and $\EE[|W|]=|A|/2$.

								Since $T$ is independent and $d_T(a)$ is odd for every $a\in W$, we can apply Lemma~\ref{lem:indep} with $P=T$ to obtain a graph $J$ with vertex set $W$ such that every even set of $J$ is contained in an odd set of $G$. Thus, by Lemma~\ref{lem:count}, we have \YZ{$e(G[W])\le120f_o^*(G)+11|W|$} for every $T\subseteq B$. Taking expectations, we obtain
									\[
									\tfrac14 e(G[A])\le\EE[e(G[W])]
									\le\YZ{120f_o^*(G)+11\EE[|W|]=120f_o^*(G)+\tfrac{11}{2}|A|,}
									\]
									and hence 
									\[\YZ{e(G[A])\le480f_o^*(G)+22|A|.}\]
								
								Choose an inclusion-minimal $A_0\subseteq A$ such that every vertex of $B$ has a neighbour in $A_0$. Let $A_1=\{a\in A\setminus A_0:d_{A_0}(a)=0\}$, and choose an inclusion-minimal $B_0\subseteq B$ such that every vertex of $A_1$ has a neighbour in $B_0$. Minimality gives private neighbours $u_a\in B$ for $a\in A_0$ and $v_b\in A_1$ for $b\in B_0$, with
								\[
								N(u_a)\cap A_0=\{a\},\qquad N(v_b)\cap B_0=\{b\}.
								\]
								Put $L_0=\{u_a:a\in A_0\}$, $P_0=\{v_b:b\in B_0\}$, $C=A_0\cup B_0$ and $D=C\cup L_0\cup P_0$.

								We now verify the hypotheses of Lemma~\ref{lem:outside} for $C$, $D$ and $V_0=V(G)\setminus D$. For each $a\in A_0$, since $B$ is independent and $u_a\in B$, we have $N(u_a)\cap C=N(u_a)\cap A_0=\{a\}$. Similarly, for each $b\in B_0$, since $v_b\in A_1$ has no neighbour in $A_0$, we have $N(v_b)\cap C=N(v_b)\cap B_0=\{b\}$. Thus every vertex of $C$ has a private neighbour in $D$.
								Every vertex outside $D$ has a neighbour in $A_0$ or belongs to $A_1$ and has a neighbour in $B_0$.
								Lemma~\ref{lem:outside} now gives 
								\[e(G[V(G)\setminus D])\le\YZ{480f_o^*(G)+22(n-|D|).}\]
								
								The edges between $A_0$ and $L_0$ form a perfect matching, as do those between $P_0$ and $B_0$.

								Every edge between $A$ and $B$ with an end in $D$ has an end in $A_0\cup L_0$ or in $P_0\cup B_0$. We apply Lemma~\ref{lem:matched} twice, first with $C=A_0$ and $L=L_0$, and then with $C=P_0$ and $L=B_0$. Since $|P_0|=|B_0|$ and $|A_0|+|B_0|=|C|$, the number of edges between $A$ and $B$ with an end in $D$ is at most
									\[
									\YZ{(484f_o^*(G)+22|A|+|A_0|)+(484f_o^*(G)+22|A|+|P_0|)
									=968f_o^*(G)+44|A|+|C|.}
									\]
								
								Since $B$ is independent, every edge of $G$ either belongs to $G[A]$ or joins $A$ and $B$. Every edge between $A$ and $B$ either has an end in $D$ or belongs to $G[V(G)\setminus D]$. Thus
								
									\[
								\YZ{	\begin{aligned}
										m&\le (480f_o^*(G)+22|A|)+(480f_o^*(G)+22(n-|D|))
										+(968f_o^*(G)+44|A|+|C|)\\
										&=1928f_o^*(G)+66|A|+22(n-|D|)+|C|\\
										&\le 1928f_o^*(G)+88n,
									\end{aligned}}
									\]
									where the last inequality follows from $|A|\le n$ and $|C|\le |D|$.
							\end{proof}
							
							\begin{proof}[Proof of Corollary~\ref{cor:main}]
								Let $n$ be the number of non-isolated vertices of $G$. As shown after Theorem~\ref{thm:zeng}, $n\le10f_o^*(G)$. By Theorem~\ref{thm:main}, $m\le \YZ{1928f_o^*(G)+880f_o^*(G)=2808f_o^*(G)}$.
							\end{proof}
							
							\section{Concluding remarks}\label{sec:remarks}
							
							We made no attempt to optimise the constants in Theorems~\ref{thm:main}, \ref{thm:even} and~\ref{thm:pairs}, and the constant in Corollary~\ref{cor:main} is certainly far from best possible. On the other hand, let $\overline{C_7}$ be the complement of the cycle $C_7$ on seven vertices, which was used by Caro~\cite{Caro} to show that the constant in the bound $f_o(G)\ge cn$ cannot be larger than $2/7$. 
							Every non-empty odd set $S$ of $\overline{C_7}$ has even size, so it is a proper subset of $V(C_7)$. At each $v\in S$ the degrees in $C_7[S]$ and $\overline{C_7}[S]$ sum to $|S|-1$. Hence $C_7[S]$ is an even forest and is therefore edgeless. An independent set of $C_7$ has at most three vertices, so $|S|=2$. Hence $f_o^*(\overline{C_7})=1$, while $\overline{C_7}$ has $14$ edges, and so the constant \YZ{$1/2808$} in Corollary~\ref{cor:main} cannot be replaced by a constant larger than $1/14$. For graphs of maximum degree at most four, $1/14$ is the correct constant. Indeed, Ai, Guo, Gutin, Hao and Yeo~\cite{AGGHY} proved that $f_o(G)\ge2n/7$ for every graph $G$ of maximum degree at most four with $n$ vertices, none of which is isolated. Such a graph has at most $2n$ edges, and an odd set with $k$ vertices induces at least $k/2$ edges, so $f_o^*(G)\ge n/7\ge e(G)/14$.
							
							\begin{question}\label{q:constant}
								What is the supremum of the constants $c$ such that $f_o^*(G)\ge c\,e(G)$ for every graph $G$? Is it equal to $1/14$?
							\end{question}
							
							For even induced subgraphs, the term $\YZ{11}n$ in Theorem~\ref{thm:even} cannot be replaced by $cn$ with a constant $c<1$, because a tree on $n$ vertices has $n-1$ edges and no even induced subgraph with an edge. It would be interesting to determine the best possible constants in Theorems~\ref{thm:even} and~\ref{thm:pairs}.
							
							\section*{Declaration of AI use} During the preparation of this manuscript, the authors used ChatGPT to assist with the exposition, to explore alternative approaches to some mathematical arguments, and to simplify and refine proofs. The authors have carefully reviewed and revised the GPT-generated suggestions and take full responsibility for the mathematical correctness and the content of the final manuscript.
							
							\section*{Acknowledgement}
							We thank Baoyindureng Wu for his excellent talk at the Nankai
							Graph Theory Conference in June 2025, which inspired us to study
							odd and even induced subgraphs and, in particular, to seek bounds
							on $f_o^*(G)$ and $f_e^*(G)$.


\begin{thebibliography}{99}
								
								\bibitem{AGGHY}
								J.~Ai, Q.~Guo, G.~Gutin, Y.~Hao and A.~Yeo.
								\newblock Odd induced subgraphs in graphs of maximum degree four.
								\newblock arXiv:2511.15489, 2025.
								
								
								\bibitem{Caro}
								Y.~Caro.
								\newblock On induced subgraphs with odd degrees.
								\newblock {\em Discrete Math.}, 132:23--28, 1994.
								
								\bibitem{CT}
								T.~M.~Cover and J.~A.~Thomas.
								\newblock {\em Elements of Information Theory}, second edition.
								\newblock Wiley, Hoboken, 2006.
								
								\bibitem{FK}
								A.~Ferber and M.~Krivelevich.
								\newblock Every graph contains a linearly sized induced subgraph with all degrees odd.
								\newblock {\em Adv. Math.}, 406:108534, 2022.
								
								\bibitem{GHZ}
								G.~Gutin, Y.~Hao and Y.~Zhou.
								\newblock Weighted counting formula and $2n/21$ lower bound for induced subgraphs with prescribed degree parities.
								\newblock arXiv:2607.16424, 2026.
								
								\bibitem{Liu}
								X.~Liu.
								\newblock Power-law lower bounds for edges in odd induced subgraphs.
								\newblock arXiv:2611.31300, 2026.
								
								\bibitem{Lovasz}
								L.~Lov\'asz.
								\newblock {\em Combinatorial Problems and Exercises}, second edition.
								\newblock North-Holland, Amsterdam, 1993.
								
								\bibitem{Scott92}
								A.~D.~Scott.
								\newblock Large induced subgraphs with all degrees odd.
								\newblock {\em Combin. Probab. Comput.}, 1:335--349, 1992.
								
								\bibitem{Scott01}
								A.~D.~Scott.
								\newblock On induced subgraphs with all degrees odd.
								\newblock {\em Graphs Combin.}, 17:539--553, 2001.
								
								\bibitem{YZ}
								X.~Yang and Q.~Zeng.
								\newblock On large odd induced subgraphs of graphs.
								\newblock {\em Ann. Appl. Math.}, 42:193--203, 2026.
								
								\bibitem{Zeng}
								Q.~Zeng.
								\newblock Large odd induced subgraphs via odd cuts.
								\newblock arXiv:2607.23266, 2026.
								
							\end{thebibliography}
						\end{document}